\documentclass[12pt]{article}

\usepackage{psfig}

\usepackage{amssymb}

\def\suffix{ps}
\newcount\system
\global\system=3   

\def\ifundefined#1{\expandafter\ifx\csname#1\endcsname\relax}
\ifundefined{figdir}\def\figdir{}\fi
\newcount\firstline
\newdimen\pswidth  \newdimen\xleft
\newdimen\psheight \newdimen\ytop \newdimen\ybot
\newcount\justx \newcount\justy
\global\justx=0 \global\justy=0
\newdimen\vpos \newtoks\labeL 
\newread\labeLfile \newdimen\xcoord \newdimen\ycoord
\newif\ifdoit 
\newbox\labox
\newdimen\xdvikwid 
\newdimen\xdvikht
\newdimen\pspoints
\newdimen\rwi
\newcount\temp
\def\readdim#1{\global\read\labeLfile to \temp
\global #1=\temp pt}
\def\figcrop#1{\par
\openin\labeLfile=\figdir#1.lbl                                              
\global\read\labeLfile to\firstline\message{#1}               
\global\read\labeLfile to\temp
\readdim{\ybot}
\readdim{\xleft}
\readdim{\ytop}
\global\read\labeLfile to\justx
\global\read\labeLfile to\justy
\global\read\labeLfile to\labeL
\readdim{\pswidth}
\global\advance\pswidth by -\xleft
\readdim{\psheight}
\global\advance\ybot by -\psheight
\global\advance\psheight by -\ytop
\global\read\labeLfile to\justx
\global\read\labeLfile to\justy
\global\read\labeLfile to\labeL
\vbox to\psheight{\vfill
\ifnum\system=1
\fi 
\ifnum\system=2
\fi 
\ifnum\system=3
                                                 \fi         
\ifnum\system=4
\fi         
\ifnum\system=1
\hbox to \pswidth{\kern-\xleft\special{postscriptfile \figdir#1.\suffix }\hfil}\fi
\ifnum\system=2
\hbox to \pswidth{\kern-\xleft\special{ps: plotfile \figdir#1.\suffix }\hfil}\fi
\ifnum\system=3
\hbox to \pswidth{\kern-\xleft\includegraphics{\figdir#1.\suffix}\hfil}\fi
\ifnum\system=4
\hbox to \pswidth{\kern-\xleft\includegraphics{\figdir#1.\suffix}\hfil}\fi
\ifnum\system=5
\hbox to \pswidth{\kern-\xleft\includegraphics{\figdir#1.\suffix}\hfil}\fi 
\ifnum\system=6
   \xdvikwid=\pswidth
   \xdvikht=\psheight
   {\global\divide\xdvikwid by \pspoints}
   {\global\divide\xdvikht by \pspoints}
   \rwi=\xdvikwid
    {\global\multiply\rwi by 10}
\hbox to \pswidth{\kern-\xleft\includegraphics{\figdir#1.\suffix\space}\hfil}\fi                   
\vskip -\baselineskip
\vskip -\ybot 
\vskip-\psheight %
\hbox to\pswidth  {\hss}%
\parindent=0pt\offinterlineskip                                       
\vpos=0 pt%
\loop\readdim{\xcoord}                                 
\ifdim \xcoord < -999pt \doitfalse\else\doittrue\fi                        
\ifdoit \advance \xcoord by -\xleft
\readdim{\ycoord}
\advance \ycoord by -\ytop                              
\global\read\labeLfile to\justx                                       
\global\read\labeLfile to\justy                                       
\global\read\labeLfile to\labeL
\global\setbox\labox=\hbox{\labeL\hskip-0.3em}%
\advance\vpos by-\ycoord                                              
\vskip-\vpos \vpos=\ycoord                                         
\hbox to\pswidth{\hskip\xcoord %
\hbox to 0pt{\ifnum\justx>0\hss\fi%
\vbox to0pt{%
\ifnum\justy<2\vss\fi%
\copy\labox\kern0pt%
\ifnum\justy>0\vss\fi}%
\ifnum\justx<2\hss\fi}%
\hss}%
\repeat%
\advance\vpos by-\psheight%
\vskip-\vpos %
}\closein\labeLfile}
\def\figplace#1#2#3{
\openin\labeLfile=\figdir#1.lbl
\ifeof \labeLfile
       \immediate\write16{***Can't find \figdir#1.lbl; Skipping it.***}
\else  \closein\labeLfile
       \null\hskip#2\raise #3 \hbox{\figcrop{#1}}
\fi
}
\def\CC{\mathbb{C}}

\newfam\gotfam
\font\twlgot=eufm10 at 12pt \font\tengot=eufm10
 \font\sevengot=eufm7
\textfont\gotfam=\twlgot \scriptfont\gotfam=\tengot
\scriptscriptfont\gotfam=\sevengot

\newtheorem{theor}   {Theorem}

\newcommand{\be}  {\begin{equation}}
\newcommand{\ee}  {\end{equation}}
\newcommand{\bea} {\begin{eqnarray}}
\newcommand{\eea} {\end{eqnarray}}
\newcommand{\lp}  {\left(}
\newcommand{\rp}  {\right)}

\newcommand{\Br}  {\overline}

\newcommand{\Om}  {\Omega}

\newcommand{\Ga}  {\Gamma}

\newcommand{\al}  {\alpha}

\newcommand{\et}  {\eta}

\newcommand{\de}  {\delta}

\newcommand{\ph}  {\phi}

\def\Br{\overline}

\newcommand{\eqdef} {\stackrel{\rm def}{=}}

\begin{document}

\title{A Physicist's Proof of the Lagrange-Good Multivariable
Inversion Formula}

\author{Abdelmalek Abdesselam \\
\\
{\small D{\'e}partement de Math{\'e}matiques}\\
{\small  Universit{\'e} Paris XIII, Villetaneuse}\\
{\small Avenue J.B. Cl{\'e}ment, F93430 Villetaneuse, France}\\
}

\maketitle


{\abstract{
We provide yet another proof of the classical Lagrange-Good
multivariable inversion formula using the techniques of quantum field
theory.
} }

 

\section{Introduction}

The Lagrange inversion formula~\cite{Lagrange} is one of the most useful
tools in enumerative combinatorics (see~\cite{GouldenJ,Stanley}).
Various efforts have been devoted to finding purely combinatorial
proofs and generalizations of this formula.
One of the many such generalizations is the extension
from the one variable to the multivariable case.
Early contributions in this direction can be found
in~\cite{Laplace,Jacobi,Darboux,Stieltjes,Poincare},
but the credit for the discovery of the general
multivariable formula is usually attributed to
the mathematical statistician I. J. Good~\cite{Good}.
We recommend~\cite{Gessel} for a clear and thorough
presentation as well as for more complete references.
Quoted from~\cite{Gessel}
the Lagrange-Good formula says the following.

\begin{theor}
Let the formal power series $f_1,\ldots, f_m$
in the variables $x_1,\ldots,x_m$ be defined by
\be
f_i=x_i g_i(f_1,\ldots,f_m)\ \ ,
\ \  
1\le i \le m
\ee
for some formal power series $g_i(x_1,\ldots,x_m)$.
Then the coefficient of
\[
\frac{x_1^{n_1}}{n_1!}\ldots
\frac{x_m^{n_m}}{n_m!}
\ \ \ 
{\rm in}
\ \ \ 
\frac{x_1^{k_1}}{k_1!}\ldots
\frac{x_m^{k_m}}{k_m!}
g_1^{n_1}(x_1,\ldots,x_m)\ldots
g_m^{n_m}(x_1,\ldots,x_m)
\]
is equal to the coefficient of
\[
\frac{x_1^{n_1}}{n_1!}\ldots
\frac{x_m^{n_m}}{n_m!}
\ \ \ 
{\rm in}
\ \ \ 
\frac{f_1^{k_1}}{k_1!}\ldots
\frac{f_m^{k_m}}{k_m!}
\times
\frac{1}{det(\de_{ij}-x_i g_{ij}(f_1,\ldots,f_m))}
\]
where
\be
g_{ij}(x_1,\ldots,x_m)
\eqdef
\frac{\partial g_i}{\partial x_j}(x_1,\ldots,x_m)\ \ .
\ee
\end{theor}

The odd-looking determinant in the denominator was probably one
of the reasons this general formula was not discovered
until~\cite{Good}. Remark however that a similar determinantal
denominator appeared earlier in the classical
MacMahon master theorem~\cite{MacMahon}.
This is no coincidence since the latter is well known to be the
linear special case of the Lagrange-Good formula.

We will use the quantum field theory model introduced
in~\cite{Abdesselam1}, which is related to an earlier
formula of G. Gallavotti~\cite{Gallavotti}
for the Lindstedt series in KAM theory,
in order to express the compositional inverse of a power
series in the multivariable setting.
Our ``proof'' of the Lagrange-Good formula will follow
from this representation of the formal inverse
by straightforward and quite natural field theoretical computations
which will, in particular, ``explain'' the
determinantal denominator as a {\em normalization factor
for a probability measure}.
We warn the mathematical reader that reckless use will be
made of wildly divergent integrals, if understood in the Lebesgue sense,
and of quantum field theory terminology.
However, when interpreted according to the formalism of our forthcoming
article~\cite{Abdesselam2}, our ``proof'' becomes a proof.
Upon closer inspection, the reader who is familiar with earlier
combinatorial proofs of the Lagrange-Good formula, like say
in~\cite{Gessel} or~\cite{EhrenborgM}, will undoubtedly have
an impression of ``d\'ej\`a vu''. Indeed, the only connected
Feynman graphs of our quantum field theory model are
either a single tree or a collection of trees branching off
a central loop.
In fact, we provided the following ``proof'', not so much for
its originality, but for its entertainment value, as another
instance of the magic of quantum field theory (see~\cite{Cartier}).

\bigskip
\noindent{\bf Aknowledgments :}
The content of this article is an outgrowth of techniques
developed in collaboration with V. Rivasseau and presented
in~\cite{Abdesselam1}.  
We thank J. Feldman for his invitation to the Mathematics Department
of the University of British Columbia where part of this work was done.
The pictures in this article have been drawn using a software package
that was kindly provided by J. Feldman.

\section{The ``proof''}
First, we avoid the use of multiindices and write
\be
g_i(x_1,\ldots,x_m)=\sum_{d\ge 0}
\frac{1}{d!}
\sum_{\al_1,\ldots,\al_d=1}^m
w_{i,\al_1\ldots\al_d}^{[d]}
x_{\al_1}\ldots x_{\al_d}
\ee
where the tensor element $w_{i,\al_1\ldots\al_d}^{[d]}$
is completely symmetric in $\al_1,\ldots,\al_d$.
Therefore the $f_i(x_1,\ldots,x_m)$ are the
solutions of
\be
f_i=x_i\sum_{d\ge 0}
\frac{1}{d!}
\sum_{\al_1,\ldots,\al_d=1}^m
w_{i,\al_1\ldots\al_d}^{[d]}
f_{\al_1}\ldots f_{\al_d}
\ee
which can be rewritten as the direct reversion problem,
with unknowns $f_1,\ldots,f_m$,
\be
y_i=\Ga_i(f)
\ee
with
$y_i\eqdef x_i w_{i}^{[0]}$ and
\be
\Ga_i(f)
\eqdef
f_i-\sum_{d\ge 1}
\sum_{\al_1,\ldots,\al_d=1}^m
\frac{1}{d!}x_i
w_{i,\al_1\ldots\al_d}^{[d]}
f_{\al_1}\ldots f_{\al_d}
\label{F1}
\ee
that is
\be
\Ga_i(f)=
\sum_{d\ge 1}
\frac{1}{d!}
\sum_{\al_1,\ldots,\al_d=1}^m
\et_{i,\al_1\ldots\al_d}^{[d]}
f_{\al_1}\ldots f_{\al_d}
\ee
with
\be
\et_{i,\al}^{[1]}\eqdef
\de_{i\al}-x_i w_{i,\al}^{[1]}
\ \ {\rm for}\ \ 
i,\al\in\{1,\ldots,m\}
\ee
and
\be
\et_{i,\al_1\ldots\al_d}^{[d]}
\eqdef
-x_i
w_{i,\al_1\ldots\al_d}^{[d]}
\ee
for $d\ge 2$, and $i,\al_1,\ldots,\al_d\in\{1,\ldots,m\}$.

It was shown in~\cite{Abdesselam1} (see~\cite{Abdesselam2} for more
detail) that the solution of such a reversion problem
is given by the perturbation expansion of the following
quantum field theory one-point function
\be
f_i=\frac{\int d{\Br \ph}d\ph\ 
\ph_i e^{-{\Br \ph}\Ga(\ph)+{\Br\ph}y}}
{\int d{\Br \ph}d\ph\ 
e^{-{\Br \ph}\Ga(\ph)+{\Br\ph}y}}\ \ .
\label{F2}
\ee
Here ${\Br\ph}_1,\ldots{\Br\ph}_m$,$\ph_1,\ldots,\ph_m$
are the components of a complex Bosonic field.
The integration is over $\CC^m$ with the measure
\be
d{\Br \ph}d\ph
\eqdef
\prod_{i=1}^m
\lp
\frac{d(Re\ \ph_i)d(Im\ \ph_i)}{\pi}
\rp
\ ,
\ee
we used the notation ${\Br \ph}\Ga(\ph)\eqdef\sum_{i=1}^m
{\Br\ph}_i \Ga_i(\ph_1,\ldots,\ph_m)$, and ${\Br\ph}y\eqdef
\sum_{i=1}^m {\Br\ph}_i y_i$.
If $\Om({\Br\ph},\ph)$ is a function of the fields, we use the notation
\be
<\Om({\Br\ph},\ph)>_U
\eqdef
\int d{\Br \ph}d\ph\ 
\Om({\Br\ph},\ph)
e^{-{\Br \ph}\Ga(\ph)+{\Br\ph}y}
\ee
for the corresponding {\em unnormalized} correlation function,
and
\be
<\Om({\Br\ph},\ph)>_N
\eqdef
\frac{1}{Z}
<\Om({\Br\ph},\ph)>_U
\ee
for the corresponding {\em normalized} correlation function,
where the normalization factor is
\be
Z\eqdef
\int d{\Br \ph}d\ph\ 
e^{-{\Br \ph}\Ga(\ph)+{\Br\ph}y}\ \ .
\ee
Finally we denote by $<\ldots>_C$
the {\em connected} correlation functions,
also known as {\em cumulants} or {\em semi-invariants}
in mathematical statistics and probability theory.

Note that the ``action'' $S({\Br\ph},\ph)\eqdef
{\Br\ph}\Ga(\ph)-{\Br\ph}y$ in the exponential
can be separtated into quadratic and nonquadratic parts by writing
\be
\Ga(\ph)=C^{-1}\ph-H(\ph)
\ee
with
\be
[C^{-1}]_{ij}\eqdef
\et_{i,j}^{[1]}=\de_{ij}-x_i w_{i,j}^{[1]}
\ee
and
\be
H(\ph)\eqdef
\sum_{d\ge 2}
\sum_{\al_1,\ldots,\al_d=1}^m
\lp -\frac{1}{d!}\rp
\et_{i,\al_1\ldots\al_d}^{[d]}
\ph_{\al_1}\ldots \ph_{\al_d}\ .
\ee
$C$ is the free propagator of our theory,
${\Br\ph}H(\ph)$ is the interaction potential
and ${\Br\ph}y$ contains the sources which can be treated
as particular vertices of the interaction.
Therefore
\be
e^{-S({\Br\ph},\ph)}
=e^{-{\Br\ph}C^{-1}\ph+{\Br\ph}H(\ph)+{\Br\ph}y}
\ee
and we let
\be
d\mu_C({\Br\ph},\ph)\eqdef
\frac{d{\Br\ph}d\ph}{det\ C}
e^{-{\Br\ph}C^{-1}\ph}
\ee
be the normalized complex Gaussian measure
with covariance $C$.
As a result
\bea
Z & = & 
\int d{\Br\ph}d\ph\ 
e^{-S({\Br\ph},\ph)} \\
 & = & (det\ C)
\int d\mu_C({\Br\ph},\ph)\ 
e^{{\Br\ph}H(\ph)+{\Br\ph}y}\ .
\eea
Now, by the standard rules of perturbative quantum field theory,
\[
\log\lp
\int d\mu_C({\Br\ph},\ph)\ 
e^{{\Br\ph}H(\ph)+{\Br\ph}y}
\rp
\]
is the sum over connected vacuum Feynman diagrams built using
the propagators
\be
\figplace{propag}{0 in}{-0.065 in}=C_{ij}
\ee
the $H$-vertices
\be
\figplace{Hvertex}{0 in}{-0.45 in} =
-\et_{i,\al_1\ldots\al_d}^{[d]}
\ee
with $d\ge 2$, and the $y$-vertices
\be
\figplace{yvertex}{0 in}{-0.47 in}=y_i\ .
\ee
These diagrams are made of a single oriented loop
of $H$-vertices linked by free propagators $C$, on which
tree diagrams terminating with $y$-vertices are hooked.
Since the sum over such tree diagrams builds the one-point
function $<\ph_i>_N=f_i=\Ga^{-1}(y)$,
it is easy to see that
\be
\log\lp
\int d\mu_C({\Br\ph},\ph)\ 
e^{{\Br\ph}H(\ph)+{\Br\ph}y}
\rp
=
\sum_{k\ge 1}
\frac{1}{k}
tr{\left[
C\partial H(\Ga^{-1}(y))
\right]}^k
\ee
where
$\partial H(z)$ is the matrix with entries
$\frac{\partial H_i}{\partial z_j}(z)$.
Therefore
\be
Z=(det\ C)
e^{-tr\log(I-C\partial H(\Ga^{-1}(y)))}
\ee
or
\bea
Z^{-1} & = & det\lp
C^{-1}(I-C\partial H(\Ga^{-1}(y)))
\rp\\
 & = & det\lp
C^{-1}-\partial H(\Ga^{-1}(y))
\rp
\eea
Now note that
\be
\partial\Ga(\ph)=C^{-1}-\partial H(\ph)
\ee
so
\be
Z^{-1}=
det\lp\partial\Ga(\Ga^{-1}(y))
\rp=
det(\partial\Ga(f))\ .
\ee
Now we also have by (\ref{F1})
\be
\Ga_i(f)=
f_i-(x_i g_i(f)-x_i w_{i}^{[0]})
\label{F3}
\ee
and thus
\bea
{[\partial\Ga(f)]}_{ij} & = &
\frac{\partial}{\partial f_j}
\lp
f_i-x_i g_i(f)+x_i w_{i}^{[0]}
\rp\\
 & = & \de_{ij}-x_i g_{ij}(f)
\eea
that is
\be
Z=
\frac{1}{det(\de_{ij}-x_i g_{ij}(f))}
\ee
which is our interpretation of the determinantal denominator
in the Lagrange-Good formula as a normalization factor for a
probability measure.

Besides, (\ref{F3}) can be rewritten as
\be
\Ga_i(f)-y_i=
f_i-x_i g_i(f)
\ee
that is (\ref{F2}) becomes
\be
f_i=\frac{1}{Z}
\int d{\Br\ph}d\ph\ \ph_i
e^{-{\Br\ph}\ph+{\Br\ph} x g(\ph)}
\ee
with
${\Br\ph} x g(\ph)\eqdef\sum_{i=1}^m
{\Br\ph}_i x_i g_i(\ph)$ and
\be
Z=\int d{\Br\ph}d\ph\ 
e^{-{\Br\ph}\ph+{\Br\ph} x g(\ph)}\ .
\ee
Now
\be
\frac{f_1^{k_1}}{k_1!}\ldots
\frac{f_m^{k_m}}{k_m!}
\times
\frac{1}{det(\de_{ij}-x_i g_{ij}(f))}
=
\frac{Z<\ph_1>_{C}^{k_1}\ldots<\ph_m>_{C}^{k_m}}
{{k_1}!\ldots{k_m}!}
\ee
but
\be
<\ph_{1}^{k_1}\ldots\ph_{m}^{k_m}>_N=
<\ph_1>_{C}^{k_1}\ldots<\ph_m>_{C}^{k_m}
\ee
because a connected graph can hook to at most one of the
sources $\ph_i$.
As a result
\be
\frac{f_1^{k_1}}{k_1!}\ldots
\frac{f_m^{k_m}}{k_m!}
\times
\frac{1}{det(\de_{ij}-x_i g_{ij}(f))}
=\int d{\Br\ph}d\ph\ 
\frac{\ph_1^{k_1}}{k_1!}\ldots
\frac{\ph_m^{k_m}}{k_m!}
e^{-{\Br\ph}\ph+{\Br\ph} x g(\ph)}\ .
\ee
Now it all becomes very simple since, on expanding $e^{{\Br\ph} x g(\ph)}$,
one gets
\bea
\lefteqn{
\frac{f_1^{k_1}}{k_1!}\ldots
\frac{f_m^{k_m}}{k_m!}
\times
\frac{1}{det(\de_{ij}-x_i g_{ij}(f))} = }
&  & \nonumber\\
 & &
\sum_{n_1,\ldots,n_m=0}^{+\infty}
\frac{x_1^{n_1}}{n_1!}\ldots
\frac{x_m^{n_m}}{n_m!}
\int d\mu_I({\Br\ph},\ph)
\prod_{a=1}^{m}
\lp
\frac{\ph_a^{k_a}}{k_a!}
\rp
\prod_{b=1}^{m}
{\lp
{\Br\ph}_b g_b(\ph)
\rp}^{n_b}
\eea
where $d\mu_I({\Br\ph},\ph)\eqdef
d{\Br\ph}d\ph\ e^{-{\Br\ph}\ph}$, the Gaussian
measure with covariance equal to the identity matrix.
Now, by integration of the ${\Br\ph}$'s by parts,
\be
\int d\mu_I({\Br\ph},\ph)
\prod_{a=1}^{m}
\lp
\frac{\ph_a^{k_a}}{k_a!}
\rp
\prod_{b=1}^{m}
{\lp
{\Br\ph}_b g_b(\ph)
\rp}^{n_b}
=
\int d\mu_I({\Br\ph},\ph)
\Om(\ph)
\ee
with
\be
\Om(\ph)
\eqdef
\frac{\partial^{n_1+\cdots+n_m}}
{\partial\ph_{1}^{n_1}\ldots\partial\ph_{m}^{n_m}}
\lp
\frac{\ph_1^{k_1}}{k_1!}\ldots
\frac{\ph_m^{k_m}}{k_m!}
g_1(\ph)^{n_1}\ldots
g_m(\ph)^{m_1}
\rp\ .
\ee
Since $\Om(\ph)$ only depends on $\ph$,
$\int d\mu_I({\Br\ph},\ph)\ 
\Om(\ph)$
is equal to the constant term of $\Om(\ph)$
which
is easily seen to be
the coefficient of
\[
\frac{\ph_1^{n_1}}{n_1!}\ldots
\frac{\ph_m^{n_m}}{n_m!}
\ \ \ 
{\rm in}
\ \ \ 
\frac{\ph_1^{k_1}}{k_1!}\ldots
\frac{\ph_m^{k_m}}{k_m!}
g_1^{n_1}(\ph)\ldots
g_m^{n_m}(\ph)
\ .
\]
which concludes our ``proof''.

\end{document}